\documentclass[final]{siamltex}

\usepackage{amsfonts,amsmath,color,verbatim}
\usepackage{stmaryrd}
\usepackage[mathscr]{eucal}
\usepackage{epsfig}
\usepackage{comment}
\usepackage[ruled,vlined]{algorithm2e}
\usepackage{color}

\def\eps{\varepsilon}
\def\e          {\ensuremath{\mathrm{e}}}

\def\C          {\ensuremath{\mathbb C}}

\def\R          {\ensuremath{\mathbb R}}

\def\Hinf{\mathcal{H}^\infty}

\newcommand{\cM}{\mathcal{M}}
\newcommand{\cS}{\mathcal{S}}

\def\iu{\mathrm{i}}

\renewcommand{\e}{{\rm e}}
\newcommand{\clambda}{\overline{\lambda}}

\def\conj{\overline}
\def\clambda{{\overline{\lambda}}}

\renewcommand{\Re}{\text{\rm Re}}
\renewcommand{\Im}{\text{\rm Im}}

\def\oeps{\eps_\star}
\def\ophi{\phi}
\def\wt{\widetilde}

\def\begcitation#1{\leftskip=2.5cm
    \noindent%
    {\small%
   \abovedisplayskip 3pt plus 1pt minus 1pt%
   \belowdisplayskip 3pt plus 1pt minus 1pt%
   #1}}
\def\signed #1 {{\unskip\nobreak\hfil\penalty50
    \hskip2em\hbox{}\nobreak\hfil\small#1
    \parfillskip=0pt \finalhyphendemerits=0 \par}}
\def\endcitation{\vskip7pt\par\leftskip=0pt} 
\def\sqr#1#2{{\vcenter{\vbox{\hrule height.#2pt
         \hbox{\vrule width.#2pt height#1pt \kern#1pt
           \vrule width.#2pt}
             \hrule height.#2pt}}}}

\newtheorem{remark}[theorem]{Remark}
\newtheorem{assumption}[theorem]{Assumption}

\title{Transient dynamics under structured perturbations: bridging unstructured and structured pseudospectra}

\author{Nicola Guglielmi\footnotemark[1] \and Christian Lubich\footnotemark[3]}
\begin{document}

\maketitle

\renewcommand{\thefootnote}{\fnsymbol{footnote}}
\footnotetext[1]{Division of Mathematics, 
Gran Sasso Science Institute,
Via Crispi 7,
I-67100   L' Aquila,  Italy. Email: {\tt nicola.guglielmi@gssi.it}}
\footnotetext[3]{Mathematisches Institut,
       Universit\"at T\"ubingen,
       Auf der Morgenstelle 10,
       D--72076 T\"ubingen,
       Germany. Email: {\tt lubich@na.uni-tuebingen.de}}
\renewcommand{\thefootnote}{\arabic{footnote}}

\begin{abstract}
The structured $\eps$-stability radius is introduced as a quantity to assess the robustness of transient bounds of solutions to linear differential equations under structured perturbations of the matrix. This applies to general linear structures such as complex or real matrices with a given sparsity pattern or with restricted range and corange, or special classes such as Toeplitz matrices. The notion conceptually combines unstructured and structured pseudospectra in a joint pseudospectrum, allowing for the use of resolvent bounds as with unstructured pseudospectra and for structured perturbations as with structured pseudospectra. We propose and study an algorithm for computing the structured $\eps$-stability radius, which solves eigenvalue optimization problems 
via suitably discretized rank-1 matrix differential equations that originate from a gradient system. The proposed algorithm has essentially the same computational cost as the known rank-1 algorithms for computing unstructured and structured stability radii. Numerical experiments illustrate the behavior of the algorithm.
\end{abstract}

\begin{keywords}
pseudospectrum, stability radius, structured $\eps$-stability radius, matrix nearness problem, eigenvalue optimization, gradient system, rank-1 dynamics.
\end{keywords}

\begin{AMS}
15A18, 65F15, 93D40
\end{AMS}

\pagestyle{myheadings} \thispagestyle{plain}
\markboth{N. GUGLIELMI,  CH. LUBICH}{Transient dynamics under structured perturbations}

\section{Introduction}
\label{sec:intro} 
\phantom{abc}\\
\vskip -1mm
\begcitation{The figure shows that
the unstructured spectral value set can be a misleading indicator of the robustness of
stability.}
\signed{D.~Hinrichsen and A.~J.~Pritchard (2005), p.\,532}
\vskip 2mm
\begcitation{By contrast, the eigenvalues that arise from structured perturbations do not bear as close a relation to the resolvent norm and may not provide much information about matrix behavior.}
\signed{L.~N.~Trefethen and M.~Embree (2005), p.\,458}
\endcitation
\noindent
\subsection{Pseudospectrum and stability radius}
For $\eps>0$, the unstructured $\eps$-pseudospectrum of a matrix $A\in \C^{n\times n}$ is defined and characterized as (see, e.g., \cite{TrE05})
\begin{align}
\Lambda_\eps(A)&= \{ \lambda \in \C \,:\, \text{$\lambda$ is an eigenvalue of $A+\Theta$ for some $\Theta\in\C^{n\times n}$ with $\| \Theta \|_F \le \eps$}\}
\nonumber
\\
&= \{ \lambda \in \C \,:\, \| (A- \lambda I)^{-1} \|_2 \ge \eps^{-1} \},
\label{ps}
\end{align}
where $\|\cdot\|_F$ is the Frobenius norm and $\|\cdot\|_2$ is the matrix 2-norm. (The matrix 2-norm might be taken also in the first line, since extremal perturbations are known to have rank 1 and the two norms are the same for rank-1 matrices.)

The pseudospectrum comes with different uses that correspond to the two lines in the above characterization:
\begin{itemize}
    \item[(P)] It provides information on changes of the spectrum under general unstructured complex perturbations;
    \item[(R)] It provides information on the norm of the resolvent. 
\end{itemize}
\pagebreak[3]
The resolvent aspect (R) leads to time-uniform bounds for the solutions of asymptotically stable linear differential equations $y'(t) = Ay(t) + f(t)$, where
all eigenvalues of $A$ have negative real part; see, e.g., \cite{HinP05,TrE05}. These bounds are reciprocal to the {\em stability radius} (or distance to instability) $\oeps=\oeps(A)$, which is the largest $\eps>0$ such that $\Lambda_\eps(A)$ has no points of positive real part. By the second line of \eqref{ps}, it follows that $\oeps^{-1}$ is the $\Hinf$-norm of the resolvent of~$A$ in the complex right half-plane:
\begin{equation}\label{stab-radius}
\frac 1\oeps = \max_{\,\Re\;\!\lambda \ge 0\,} \| (A-\lambda I)^{-1} \|_2.
\end{equation}

\subsection{Structured pseudospectrum and structured stability radius}
If the emphasis is put on the perturbation aspect (P), the concept needs to be adjusted when only structured perturbations are of interest, e.g., real perturbations or perturbations with a given sparsity pattern, or perturbations with restricted range and corange, or Toeplitz matrices. For a given structure space $\cS\subset \C^{n\times n}$, which may be an arbitrary complex-linear or real-linear subspace of $\C^{n\times n}$, the {\em $\cS$-structured $\delta$-pseudospectrum} is defined in the following way (cf., e.g., \cite{HinP05}), where only structured perturbations $\Delta \in \cS$ are allowed:
\begin{equation}\label{s-ps}
\Lambda_\delta^\cS(A)= \{ \lambda \in \C \,:\, \text{$\lambda$ is an eigenvalue of $A+\Delta$ for some $\Delta\in\cS$ with $\| \Delta \|_F \le \delta$}\}.
\end{equation}
(Here it makes a difference whether the Frobenius norm or the matrix 2-norm is chosen. For the following we prefer to work with the Frobenius norm, since it is an inner-product norm.)

In contrast to the complex unstructured pseudospectrum, no characterization in terms of resolvent norms is available for structured pseudospectra.
The {\em $\cS$-structured stability radius} $\delta_\star^\cS$ is the largest $\delta>0$ such that $\Lambda_\delta^\cS(A)$ has no points of positive real part. This radius still indicates up to which perturbation size a perturbed linear dynamical system is guaranteed to remain asymptotically stable as $t\to\infty$, but it provides no information on transient solution bounds.

\subsection{Joint unstructured/structured pseudospectrum and structured $\eps$-stability radius}
\label{subsec:s-eps-sr}
We will derive transient bounds that are robust under structured perturbations. They are obtained by combining both unstructured and structured pseudospectra in a joint pseudospectrum:
\begin{align}
\nonumber
\Lambda_{\delta,\eps}^{\cS}(A) &= \{ \lambda \in \C \,:\, \text{$\lambda\in \Lambda_\eps(A+\Delta)$ for some $\Delta\in\cS$ with $\| \Delta \|_F \le \delta$} \}
\\
\label{uss-ps}
&= \{ \lambda \in \C \,:\, \text{$\lambda$ is an eigenvalue of $A+\Delta+ \Theta$ for some} 
\\[-1mm]
&\hskip 2cm \text{$\Delta\in\cS$ with $\| \Delta \|_F \le \delta$ and $\Theta \in \C^{n\times n}$ with $\| \Theta \|_F \le \eps$} 
\}
\nonumber
\\
&=  \{ \lambda \in \C \,:\, \text{$\| (A+\Delta-\lambda I)^{-1} \|_2 \ge \eps^{-1}$ for some $\Delta\in\cS$ with $\| \Delta \|_F \le \delta$} \},
\nonumber
\end{align}
where the structure space $\cS$ is again a complex-linear or real-linear subspace of $\C^{n\times n}$.
The equalities follow from the characterization in \eqref{ps}. Note that $\Lambda_{\delta,0}^{\cS}(A)=\Lambda_\delta^\cS(A)$ and $\Lambda_{0,\eps}^{\cS}(A)=\Lambda_\eps(A)$  and therefore 
$ \Lambda_\delta^\cS(A)\cup \Lambda_\eps(A)  \subset \Lambda_{\delta,\eps}^{\cS}(A)$.

A basic notion in this context is the structured {$\eps$-stability} radius defined as follows.
Here, the matrix $A\in \C^{n\times n}$ is assumed to have all eigenvalues of negative real part, and $\oeps>0$ is its stability radius.

\begin{definition} \label{def:s-eps-sr}
For $0<\eps<\oeps$, the {\em $\cS$-structured $\eps$-stability radius} of $A$, 
denoted $\delta_\eps^\cS=\delta_\eps^\cS(A)$, is the largest $\delta>0$ such that for every structured perturbation $\Delta\in\cS$ with $\| \Delta \|_F \le \delta$, the unstructured $\eps$-pseudospectrum $\Lambda_\eps(A+\Delta)$ has no points with positive real part. 
\end{definition}

The $\cS$-structured $\eps$-stability radius of $A$ 
is thus the largest $\delta>0$ such that $\Lambda_{\delta,\eps}^{\cS}(A)$ has no points with positive real part.
Note that for $\eps\to 0$, the structured $\eps$-stability radius $\delta_\eps^\cS$ 
becomes the structured stability radius $\delta_\star^\cS$. The unstructured $\eps$-stability radius is simply $\oeps-\eps$ but depending on the structure $\cS$, the $\cS$-structured $\eps$-stability radius can be significantly larger.

By the second line of \eqref{uss-ps}, the structured $\eps$-stability radius is characterized as the largest $\delta>0$ such that all eigenvalues of $A+\Delta + \Theta$ have nonpositive real part for every $\Delta\in \cS$ with $\|\Delta\|_F \le \delta$
and every $\Theta\in \C^{n\times n}$ with $\|\Theta\|_F \le \eps$. This characterization will be used in the algorithm that we propose for computing the structured $\eps$-stability radius.

On the other hand, by the last line of \eqref{uss-ps} it follows that with $\delta=\delta_\eps^\cS$,
\begin{equation}
    \label{robust-res-bound}
\frac 1\eps = \max_{\Delta\in\cS, \| \Delta\|_F \le \delta}\max_{\,\Re\lambda \ge 0\,} \| (A+\Delta-\lambda I)^{-1} \|_2,   
\end{equation}
so that the resolvent norm of $A+\Delta$  is bounded in the right complex half-plane by $\eps^{-1}$ uniformly for all
$\Delta\in\cS$ with $\|\Delta\|_F \le \delta$.
This characterization allows us to obtain robust transient bounds reciprocal to $\eps$ for linear dynamical systems with perturbed matrices $A+\Delta$ with $\Delta\in \cS$ and $\| \Delta \|_F \le \delta_\eps^\cS$. 

In \eqref{robust-res-bound} we considered $\delta=\delta_\eps^\cS$ as a function of $\eps$. Conversely, for a given $\delta$ with $0<\delta<\delta_\star^\cS$, it is of interest to know for which $\eps$ equation \eqref{robust-res-bound} holds true. We thus ask two basic questions:
\begin{itemize}
    \item Up to which size of structured perturbations are the resolvent norms of the perturbed matrices within a given bound in the right complex half-plane?
    \item For a given size of structured perturbations, what is the smallest common bound for the resolvent norms of the perturbed matrices in the right complex half-plane?
\end{itemize}
{\em In this paper we propose and study an algorithm for computing the structured $\eps$-stability radius $\delta_\eps^\cS$.}  This algorithm turns out to have essentially the same computational cost as the  algorithm for computing the stability radius based on \cite{GL11} (see also \cite{GO11} and 
\cite{GLS23}), which is known to be particularly efficient for large sparse matrices. The algorithm requires only a very minor modification to compute the resolvent bound $\eps^{-1}$ that corresponds to a given $\delta>0$ via \eqref{robust-res-bound}. So essentially the same algorithm addresses both questions.

\subsection{Outline}

In section 2 we give time-uniform bounds for solutions of homogeneous and inhomogeneous linear differential equations. These bounds are robust under structured perturbations of norm up to the structured $\eps$-stability radius. They are based on the robust resolvent bound \eqref{robust-res-bound}.

In section 3 we describe a two-level approach to computing the structured $\eps$-stability radius. The inner iteration requires solving an eigenvalue optimization problem for a pair of unstructured and structured perturbations of fixed norms, which maximize the real part of eigenvalues for given perturbation sizes $\eps$ and~$\delta$. We study a norm-constrained gradient flow for this problem and 
find that  
the unstructured component of the minimizer is of rank 1 and the structured component is a real multiple of the orthogonal projection of this rank-1 matrix onto the structure space $\cS$. In sections 4 and 5 we make use of this form of minimizers to reduce the norm-constrained gradient flow system from $\C^{n\times n}\times \cS$ to complex $n\times n$ matrices of rank~1. This translates into a system of differential equations for two $n$-vectors, which coincide with left and right eigenvectors of the extremally perturbed matrix in stationary points.
This system of differential equations is then solved numerically into a stationary point, using a splitting method that is described in section 6.


In section 7 we briefly discuss the outer iteration. This requires computing the zero of a univariate nonlinear function, which is given by the real part of rightmost eigenvalues for extremal perturbations obtained from the inner iteration, now considered as a function of the perturbation size $\delta$ or $\eps$. 
We use a Newton-type method for which we obtain a simple expression for the derivative of this function. 

In section 8 we present numerical experiments that illustrate the behavior of the proposed algorithm for computing the structured $\eps$-stability radius. We present illustrative numerical results for a small banded Toeplitz matrix and results of numerical experiments with large sparse matrices from the Matrix Market \cite{BPR97}.

\section{Transient bounds that are robust under structured perturbations}

We give two results of robust bounds of solutions of linear differential equations that follow from the robust resolvent bound \eqref{robust-res-bound}. Although the arguments used in the proofs of these bounds are not new, we include the short proofs for the convenience of the reader. The first result deals with homogeneous linear differential equations, the second result with inhomogeneous linear differential equations with zero initial value. 
In both results,
$\cS\subset \C^{n\times n}$ is a complex-linear or real-linear structure space, $A\in \C^{n\times n}$ is a given matrix with all eigenvalues of negative real part, and $0<\eps<\oeps$ where $\oeps$ is the (unstructured) stability radius of $A$ as in \eqref{stab-radius}.
Furthermore, $\delta_\eps^\cS(A)$ is the $\cS$-structured $\eps$-stability radius of $A$ as introduced in Definition~\ref{def:s-eps-sr}.

The first result is a variant of \cite[Theorem 15.2]{TrE05} with structured perturbations.
\begin{proposition} \label{lem:exp-bound-S}
For every $\Delta\in\cS$ with $\|\Delta\|_F\le \delta_\eps^\cS(A)$,
$$
\bigl\| \e^{t(A+\Delta)} \bigr\|_2  \le  
\frac{|\Gamma|}{2\pi\eps} \qquad\text{for all }\ t>0,
$$
where $\Gamma$ is a closed contour in the closed left complex half-plane that is a union of (i) the part in the left half-plane of a contour (or union of several contours) that surrounds the pseudospectrum $\Lambda_{\eps+\delta}(A)$ with $\delta=\delta_\eps^\cS(A)$ and (ii) one or several intervals on the imaginary axis that close the contour. Moreover, $|\Gamma|$ is the length of $\Gamma$.
\end{proposition}

\begin{proof} Let $\Delta\in\cS$ with $\|\Delta\|_F\le \delta= \delta_\eps^\cS(A)$.
The bound follows from the Cauchy integral representation
$$
\e^{t(A+\Delta)} = \frac1{2\pi\iu} \int_{\Gamma} \e^{t \lambda}\, (\lambda I -A-\Delta)^{-1}\, d\lambda
$$
on noting that $|\e^{t \lambda}|\le 1$ for all $\lambda\in\Gamma$ and
\begin{equation}\label{res-bound-Gamma}
\| (\lambda I -A-\Delta)^{-1} \|_2  \le  \frac1\eps \quad\ \text{ for all } \lambda\in\Gamma.
\end{equation}
This resolvent bound holds true because 
(i) the inclusion $\Lambda_{\eps}(A+\Delta) \subset \Lambda_{\eps+\delta}(A)$ implies the bound \eqref{res-bound-Gamma}
for $\lambda$ in the closure of $\C \setminus \Lambda_{\eps+\delta}(A)$, and (ii) the resolvent bound \eqref{robust-res-bound} implies the bound \eqref{res-bound-Gamma} for all $\lambda$ on the imaginary axis.
\end{proof}

The next result can be viewed, in the spirit of systems and control theory (see, e.g., \cite{HinP05}), as a bound of the input-output relation for perturbed transfer functions $H_\Delta(\lambda)= (\lambda I - A - \Delta)^{-1}$ with structured perturbations $\Delta\in\cS$.
\begin{proposition}\label{lem:ode-inhom-bound-S}
For all perturbations $\Delta\in\cS$ with $\|\Delta\|_F\le \delta_\eps^\cS(A)$,
solutions to the inhomogeneous linear differential equations 
$$
y'_\Delta(t)=(A+\Delta)y_\Delta(t)+ f(t), \qquad y_\Delta(0)=0, 
$$
share the bound (with $\|\cdot\|$ the Euclidean norm on $\C^n$)
\begin{equation} \label{robust-L2-bound}
  \biggl(  \int_0^T \| y_\Delta(t) \|^2 \, dt \biggr)^{1/2} \le \frac1\eps \,
   \biggl(  \int_0^T \| f(t) \|^2 \, dt \biggr)^{1/2} , \qquad 0\le T \le \infty.
\end{equation}
\end{proposition}

\begin{proof} We extend  $y_\Delta(t)$ and $f(t)$ to $t<0$ by zero. Their Fourier transforms $\widehat y_\Delta$ and $\widehat f$ are then related by $\iu\omega\,\widehat y_\Delta(\omega)=(A+\Delta)\widehat y_\Delta(\omega) + \widehat f(\omega)$ for all $\omega\in\R$, i.e.,
$$
\widehat y_\Delta (\omega)=(\iu\omega I - A-\Delta)^{-1}\widehat f(\omega), \qquad \omega \in \R.
$$
Applying the Plancherel formula twice yields
\begin{align*}
    &\int_\R \| y_\Delta(t) \| ^2 \, dt = \int_\R \| \widehat y_\Delta (\omega)\| ^2\, d\omega = 
    \int_\R \| (\iu\omega I - A-\Delta)^{-1}\widehat f(\omega) \| ^2 \, d\omega
    \\
    &\le \max_{\omega\in\R} \| (\iu\omega I - A-\Delta)^{-1} \|_2^2 \, \int_\R \| \widehat f (\omega)\| ^2\, d\omega =
    \max_{\omega\in\R} \| (\iu\omega I - A-\Delta)^{-1} \|_2^2 \,  \int_\R \| f(t) \| ^2 \, dt.
\end{align*}
Using \eqref{robust-res-bound} and the causality property that $y_\Delta(t)$, for $0\le t \le T$, only depends on $f(\tau)$ with $0\le \tau \le t \le T$ (this allows us to extend $f(t)$ by $0$  for $t>T$), we obtain the bound
\eqref{robust-L2-bound}.
\end{proof}

\section{Eigenvalue optimization problem and constrained gradient flow}
\label{sec:opt}
Our numerical approach to computing the structured $\eps$-stability radius $\delta_\eps^\cS(A)$ uses a two-level iteration, which solves an eigenvalue optimization problem in the inner iteration and uses a one-dimensional root-finding procedure for the outer iteration. In this section we first describe the two-level approach and then study a structure- and norm-constrained gradient flow for the eigenvalue optimization problem. This gives us useful insight into properties of the solutions of the optimization problem that will allow us to derive computationally more efficient approaches in later sections, where the norm-constrained gradient flow on $\C^{n\times n}\times \cS$ is ultimately reduced to a norm-constrained gradient flow on the manifold of $n\times n$ matrices of rank 1, which is equivalent to a system of differential equations for two vectors in $\C^n$ of unit norm.

\subsection{Two-level approach}
For any square matrix $M$, let  $\lambda(M)$ be an eigenvalue of $M$ of maximal real part (and in case there are several such eigenvalues, take, e.g., the one with maximal imaginary part).
For $\eps>0$ and $\delta>0$ we introduce the functional
\begin{equation}\label{F-eps-delta}
F_{\eps,\delta}(E,E^\cS) = - \Re \,\lambda(A+\eps E + \delta E^\cS)
\end{equation}
for $E \in \C^{n\times n}$ and $E^\cS\in \cS$, both of unit Frobenius norm. With this functional we follow a two-level approach:
\begin{itemize}
\item {\bf Inner iteration} (eigenvalue optimization): For a given $\eps>0$ and $\delta>0$, we aim to compute  matrices $E_{\eps,\delta} \in \C^{n\times n}$ and $E^\cS_{\eps,\delta} \in\cS$, both of unit Frobenius norm,
that minimize $F_{\eps,\delta}$:
\begin{equation} \label{E-ES-min}
(E_{\eps,\delta},E^\cS_{\eps,\delta}) = \arg\min\limits_{\genfrac{}{}{0pt}{}{E\in \C^{n\times n}, E^\cS \in \cS}{\| E \|_F = \| E^\cS\|_F  = 1}}
F_{\eps,\delta}(E,E^\cS).
\end{equation}

\item {\bf Outer iteration} (root finding): We compute $\delta_\eps>0$ as the smallest positive zero of the univariate function
$\ophi_\eps(\delta)=  F_{\eps,\delta}(E_{\eps,\delta},E^\cS_{\eps,\delta})$:
\begin{equation} \label{zero-delta}
\ophi_\eps(\delta_\eps)= 0.
\end{equation}
\end{itemize}
Provided that these computations succeed, we have that 
$\delta_\eps$ is the $\cS$-structured $\eps$-stability radius of $A$. Hence, $\delta=\delta_\eps$ satisfies \eqref{robust-res-bound} for the given $\eps$.

For the second question at the end of section~\ref{subsec:s-eps-sr},  the  outer iteration 
computes instead $\eps_\delta>0$ as the smallest positive zero of the function
$\psi_\delta(\eps)=  F_{\eps,\delta}(E_{\eps,\delta},E^\cS_{\eps,\delta})$:
\begin{equation} \label{zero-eps}
\psi_\delta(\eps_\delta)= 0.
\end{equation}
Then, $\eps=\eps_\delta$ satisfies \eqref{robust-res-bound} for the given $\delta$.

\subsection{Orthogonal projection onto the structure}\label{subsec:proj-structure}
Given two complex $n \times n$ matrices, we denote by  
\begin{equation*}
\langle X,Y \rangle =\sum_{i,j} \conj{x}_{ij}y_{ij} = {\rm trace}(X^* Y)
\end{equation*} 
the inner product in $\C^{n\times n}$ that induces the Frobenius norm $\| X \|_F = \langle X,X \rangle^{1/2}$.

Let $\Pi^\cS:\C^{n\times n}\to\cS$ be the orthogonal projection 
onto~$\cS$\/: for every $Z\in \C^{n\times n}$, 
\begin{equation}\label{Pi-S}
\Pi^\cS Z \in \cS \quad\text{ and } \quad \Re\langle \Pi^\cS Z, W \rangle = \Re\langle Z,W \rangle \quad \text{for all } \, W\!\in\cS.
\end{equation}
For a complex-linear subspace $\cS$, taking the real part of the complex inner product can be omitted (because with $W\in\cS$, then also $\iu W\in\cS$), but taking the real part is needed for real-linear subspaces. Note that for $\cS=\R^{n\times n}$,  we have $\Pi^\cS Z=\Re\, Z$ for all $Z\in\C^{n\times n}$.

If $\cS$ is the space of complex matrices with a prescribed sparsity pattern, then $\Pi^\cS Z$ leaves the entries of $Z$ on the sparsity pattern unchanged and annihilates those outside the sparsity pattern.  If $\cS$ is the space of real matrices with a prescribed sparsity pattern, then $\Pi^\cS Z$ takes instead the real part of the entries of $Z$ on the sparsity pattern.
We further refer to \cite[Example~2.2]{GLS23} for the projection in the case where the structure $\cS$ consists of matrices with fixed range and corange.

\subsection{Gradient of the functional}
The following result is proved in the same way as Lemma~2.4 and equation (2.10) in \cite{GLS23}, based on the derivative formula for simple eigenvalues as given, e.g., in \cite[Theorem~1]{GreLO20}.
\begin{lemma}[Gradient] \label{lem:gradient}
Let $(E(t),E^\cS(t))\in \C^{n\times n}\times \cS$, for $t$ near $t_0$, be a continuously differentiable path of matrices, with the derivative denoted by $(\dot E(t),\dot E^\cS(t))$.
Assume that $\lambda(t)$ is a simple eigenvalue of  $A+\eps E(t)+ \delta E^\cS(t)$ depending continuously on $t$, with associated left and right eigenvectors
$x(t)$ and $y(t)$ normalized such that they are of unit norm and with real and positive inner product. Let the eigenvalue condition number be
\begin{equation*}
\kappa(t) = \frac1{x(t)^* y(t)} > 0.
\end{equation*}
Then, $F_{\eps,\delta}(E(t),E^\cS(t))=-\Re\,\lambda(t)$ 
is continuously differentiable w.r.t. $t$ and, omitting the ubiquitous  argument $t$ on the righthand side,
\begin{align}
 \label{deriv}
\frac1{\kappa(t) } \,\frac{d}{dt} F_{\eps,\delta}(E(t),E^\cS(t)) &= \Re \,\bigl\langle  G_{\eps,\delta}(E,E^\cS),  \eps\dot E \bigr\rangle \\
\nonumber
&\quad + \Re\, \bigl\langle  \Pi^\cS G_{\eps,\delta}(E,E^\cS), \delta \dot E^\cS \bigr\rangle 
\end{align}
with the rank-1 matrix
\begin{equation} \label{freegrad}
G_{\eps,\delta}(E,E^\cS) = - x y^* \in \C^{n\times n} .
\end{equation}
\end{lemma}


\subsection{Norm-constrained gradient flow}
We consider the system of differential equations,
with $G=G_{\eps,\delta}(E,E^\cS)$ for short,
\begin{equation} \label{ode-E-ES}
    \begin{aligned}
     \eps \dot E &= -  G + \Re\langle G, E \rangle E,
       \\[1mm]
   \delta \dot E^\cS &= - \Pi^\cS G + \Re\langle \Pi^\cS G, E^\cS \rangle E^\cS.
    \end{aligned}
\end{equation}
Along its solutions, the unit Frobenius norms of $E(t)\in \C^{n\times n}$ and $E^\cS(t)\in\cS$  are preserved: 
$$
\frac\eps 2 \frac d{dt} \| E(t) \|_F^2 = \Re \langle \eps\dot E, E \rangle = -\Re \langle G,E \rangle 
+ \Re \langle G,E \rangle \, \|E\|_F^2 =0 ,
$$
and hence $\| E(t) \|_F^2=1$ for all $t$. By the same argument, also $\| E^\cS(t) \|_F^2=1$ for all~$t$.

The following result shows that the functional $F_{\eps,\delta}$ decays along solutions of \eqref{ode-E-ES}.

\begin{theorem}[Monotone decay of the functional]
\label{thm:monotone}
Let $(E(t), E^\cS(t))\in \C^{n\times n}\times \cS$ with $\| E(t)\|_F=\| E^\cS(t)\|_F= 1$ be a solution of the differential equation {\rm (\ref{ode-E-ES})}.
Assume that the rightmost eigenvalue $\lambda(A+\eps E(t)+\delta E^\cS(t))$ is simple and depends continuously on $t$. 
Then,
\begin{equation}
\frac{d}{dt} F_{\eps,\delta} (E(t),E^\cS(t))  \le  0,
\label{monotone}
\end{equation}
and the inequality is strict unless $(E(t),E^\cS(t))$ is a stationary point of \eqref{ode-E-ES}.
\end{theorem}

\begin{proof}
We abbreviate $G=G_{\eps,\delta}(E,E^\cS)$  and take the inner product of \eqref{ode-E-ES} with $(\dot E,\dot E^\cS)$. Using that
$\Re\langle E, \dot E \rangle =0$ and $\Re\langle E^\cS, \dot E^\cS \rangle =0$, we find
\begin{align*}
 \| \eps \dot E \|_F^2 &=  \Re \langle -G + \Re\langle G,E \rangle E, \eps \dot E\rangle =
- \Re \langle G, \eps \dot E\rangle
\\ 
\| \delta \dot E^\cS \|_F^2 &=  \Re \langle -\Pi^\cS G + \Re\langle \Pi^\cS G,E^\cS \rangle E^\cS, \delta\dot E^\cS\rangle =
- \Re \langle \Pi^\cS G, \delta\dot E^\cS\rangle
\end{align*}
and hence Lemma~\ref{lem:gradient} and \eqref{ode-E-ES} yield
\begin{align}\nonumber
\frac1{\kappa(t)} \,\frac{d}{dt} F_{\eps,\delta}(E(t),E^\cS(t)) &= \Re \langle G, \eps\dot E\rangle + \Re \langle \Pi^\cS G, \delta \dot E^\cS\rangle
\\
&= - \| \eps\dot E \|_F^2 - \| \delta\dot E^\cS \|_F^2 \le 0,
\label{c-s}
\end{align}
which on inserting \eqref{ode-E-ES} shows the decay rate of $F_{\eps,\delta}$ along a trajectory.
\end{proof}

The stationary points of the differential equation  \eqref{ode-E-ES} are characterized as follows.

\begin{theorem}[Stationary points] 
\label{thm:stat}
Let $(E_\star,E_\star^\cS)\in \C^{n\times n}\times \cS$ with $\| E_\star\|_F=\| E_\star^\cS\|_F= 1$ 
and set $G_\star=G_{\eps,\delta}(E_\star,E_\star^\cS)=-x_\star y_\star^*$.
Assume that 
\begin{itemize}
\item[(i)]   The rightmost eigenvalue $\lambda(A+\eps E+\delta E^\cS)$ is simple at $(E_\star,E_\star^\cS)$ and depends continuously on $(E,E^\cS)$ in a neighborhood of $(E_\star,E_\star^\cS)$.
\item[(ii)]  $\Pi^\cS G_\star \ne 0$.
\end{itemize}
Let $(E(t), E^\cS(t))\in \C^{n\times n}\times \cS$ be the solution of \eqref{ode-E-ES} passing through $(E_\star,E_\star^\cS)$.
Then the following statements are equivalent: 
\begin{itemize}
\item[1.] $\displaystyle\frac{ d }{dt} F_{\eps,\delta}(E(t),E^\cS(t))  = 0$. \\
\item[2.] $\dot E = 0$ and $\dot E^\cS = 0$. \\[-1mm]
\item[3.] $E_\star$ is a real multiple of $G_\star$, and $E_\star^\cS$ is a real multiple of $\Pi^\cS G_\star$. 
\end{itemize}
\end{theorem}
\medskip

\begin{proof} By \eqref{ode-E-ES}, 3.~implies 2., which implies 1. Furthermore, \eqref{c-s} and \eqref{ode-E-ES} show that 1. implies 3.
\end{proof}

In the nonsmooth nonconvex optimization problem \eqref{E-ES-min}, there can be several stationary points, and not each of them is a global minimum. Computing several trajectories with different starting values reduces the risk of getting caught in a nonoptimal local minimum. Our numerical experience indicates, however, that the global minimum can usually be found by a single discrete trajectory.

We note the following important consequence of Theorem~\ref{thm:stat}. 
\begin{corollary} \label{cor:min}
Minimizers $(E_\star,E_\star^\cS)\in \C^{n\times n}\times \cS$ with $\| E_\star\|_F=\| E_\star^\cS\|_F= 1$  of the eigenvalue optimization problem \eqref{F-eps-delta}--\eqref{E-ES-min} have the following properties under the nondegeneracy conditions (i) and (ii) of Theorem~\ref{thm:stat}:
$$
\text{
$E_\star$ is of rank~1 \quad and \quad $E_\star^\cS=\pm \eta\,\Pi^\cS E_\star$ with $\eta=1/\|\Pi^\cS E_\star\|_F$.
}
$$
\end{corollary}

This motivates us to restrict the search for a minimum of \eqref{F-eps-delta} to the reduced set of pairs of matrices $(E,E^\cS)\in \C^{n\times n}\times \cS$ with $\| E\|_F=\| E^\cS\|_F= 1$ for which we have $E^\cS=\pm \eta\,\Pi^\cS E$ with $\eta=1/\|\Pi^\cS E\|_F$ (section~\ref{sec:reduced-opt})
and for which $E$ has rank 1 (section~\ref{sec:rank-one}).

\section{Reduced optimization problem}
In view of Corollary~\ref{cor:min} we consider, instead of \eqref{F-eps-delta}, the minimization of the restricted functional
\label{sec:reduced-opt}
\begin{equation} \label{eq:redF}
\wt F_{\eps,\delta}(E) = - \Re \,\lambda\biggl(A+\eps E \pm  \delta \frac{\Pi^\cS E}{\| \Pi^\cS E\|_F}\biggr),
\end{equation}
for which we follow the program of the previous section.
For ease of presentation we consider here only the $+$ case. The $-$ case is completely analogous (just replace $\delta$ by $-\delta$ on every occurrence).
\begin{lemma}[Reduced gradient] 
\label{lem:red-gradient}
Let $E(t)\in \C^{n, n}$, for real $t$ near $t_0$, be a continuously differentiable path of matrices, with the derivative denoted by $\dot E(t)$.
Assume that $\lambda(t)$ is a simple eigenvalue of  $A+\eps E(t)+\delta \eta(t) \Pi^\cS E(t)$ depending continuously on~$t$, with
$\eta(t)=1/\| \Pi^\cS E(t)\|_F$,  with associated left and right eigenvectors
$x(t)$ and $y(t)$ of unit norm and with positive inner product. Let the eigenvalue condition number be 
$\kappa(t) = 1/(x(t)^* y(t)) > 0$.
%
Then, $\wt F_{\eps,\delta}(E(t))= \Re \, \lambda(t)$ 
is continuously differentiable w.r.t. $t$ and we have
\begin{equation} \label{deriv-red}
\frac1{ \kappa(t) } \,\frac{d}{dt} \wt F_{\eps,\delta}(E(t)) = \Re \,\bigl\langle  \wt G_{\eps,\delta}(E(t)),  \dot E(t) \bigr\rangle,
\end{equation}
where the (rescaled) gradient of $\widetilde F_{\eps,\delta}$ is the matrix, with $G=- xy^*$ and $\eta=1/\| \Pi^\cS E\|_F$,
\begin{equation} \label{red-grad}
\wt G_{\eps,\delta}(E) =  \eps G + \delta\eta\, \Pi^\cS G - \delta \eta \,\Re\bigl\langle G,\eta\,\Pi^\cS E\bigr\rangle \,\eta\,\Pi^\cS E \in \C^{n, n}.
\end{equation} 
\end{lemma}

\begin{proof}
By the derivative formula of simple eigenvalues \cite[Theorem~1]{GreLO20},
\begin{align*} 
\frac{ d }{dt} \wt F_{\eps,\delta} \left( E(t) \right) = - \Re\,\dot \lambda(t) =
-\frac{ 1}{x^*y} \,\Re \left( x^* \bigl(\eps\dot{E}+ \delta \dot \eta\, \Pi^\cS E + \delta \eta\, \Pi^\cS \dot E\bigr) y \right),
\end{align*}
where  we omit the dependence on $t$ on the right-hand side.
Since
$\Re ( x^* Z y ) = \Re\, \langle x y^* , Z  \rangle$
for every matrix $Z\in \C^{n\times n}$ and
$$
\dot \eta = \frac d{dt} \bigl\langle \Pi^\cS E , \Pi^\cS E \bigr\rangle^{-1/2} =
-\eta^3 \, \Re\,\bigl\langle \Pi^\cS E , \Pi^\cS \dot E \bigr\rangle,
$$
we obtain 
\begin{align*}
- \frac1\kappa \,\Re\,\dot \lambda & = 
\Re\, \Bigl\langle G, \eps \dot E + \delta \eta\, \Pi^\cS \dot E 
- \delta \eta^3 \,\Re\bigl\langle \Pi^\cS E, \Pi^\cS \dot E \bigr\rangle\, \Pi^\cS E \Bigr\rangle
\\
&= \Re\, \bigl\langle G, \eps \dot E + \delta \eta\, \Pi^\cS \dot E \bigr\rangle
- \delta \eta^3 \,\Re\bigl\langle \Pi^\cS E, \Pi^\cS \dot E \bigr\rangle \, 
\Re \bigl\langle G, \Pi^\cS E \bigr\rangle
\\
&= \Re\, \Bigl\langle \eps G + \delta \eta\, \Pi^\cS G 
- \delta\eta^3 \,\Re \bigl\langle G, \Pi^\cS E \bigr\rangle \,\Pi^\cS E, \dot E \Bigr\rangle,
\end{align*}
which is the stated result.
\end{proof}

We consider the norm-constrained gradient flow
\begin{equation}\label{ode-E}
\dot E = -\wt G_{\eps,\delta}(E) + \Re \,\langle \wt G_{\eps,\delta}(E), E \rangle\, E.
\end{equation}
By construction of this differential equation, we have 
$$\frac12\,\frac{d}{dt} \| E(t)\|_F^2 = \Re \langle E(t), \dot E(t) \rangle =0
$$ 
along its solutions, and so the Frobenius norm $1$ is conserved.
Since we follow the admissible direction of steepest descent of the functional $\wt F_{\eps,\delta}$ along solutions $E(t)$ of this differential equation, we have the following monotonicity property.

\begin{theorem}[Monotone decay of the restricted functional] 
\label{thm:monotone-red}
Let $E(t)\in \C^{n\times n}$ of unit Frobenius norm satisfy the differential equation {\rm (\ref{ode-E})}.

Assume that the rightmost eigenvalue $\lambda(A+\eps E(t)+\delta\eta(t)\,\Pi^\cS E(t))$ with the normalizing factor $\eta(t) = 1/\| \Pi^\cS E(t) \|_F$ is simple and depends continuously on $t$.
Then, we have
\begin{equation} 
\frac{d}{dt} \wt F_{\eps,\delta} (E(t))  \le  0.
\label{monotone-red}
\end{equation}
\end{theorem}

\begin{proof} We abbreviate $\wt G=\wt G_{\eps,\delta}(E)$  and take the inner product of \eqref{ode-E} with $\dot E$. Since
$\Re\langle E, \dot E \rangle =0$, we find
$$
\| \dot E \|_F^2 = - \Re \big\langle \wt G - \Re\langle \wt G,E \rangle E, \dot E \big\rangle =
- \Re \langle \wt G, \dot E\rangle,
$$
and hence Lemma~\ref{lem:red-gradient} and \eqref{ode-E} yield
\begin{equation}\label{c-s-red}
\frac1{\kappa} \,\frac{d}{dt} F_\eps(E(t)) = \Re \langle \wt G, \dot E\rangle = - \| \dot E \|_F^2 =
- \| \wt G - \Re\,\langle \wt G, E \rangle E \|_F^2 \le 0,
\end{equation}
which is the result as stated.
\end{proof}

The stationary points of the differential equation  \eqref{ode-E} are characterized as follows.

\begin{theorem}[Stationary points] 
\label{thm:stat-red}
Let $E_\star\in\C^{n\times n}$ with $\| E_\star\|_F=1$ be such that 
\begin{itemize}
\item[(i)]   the rightmost eigenvalue $\lambda(A+\eps E+ \delta\eta\,\Pi^\cS E)$ with $\eta = 1/\| \Pi^\cS E \|_F$ is simple at $E_\star$ and depends continuously on $E$ in a neighborhood of $E_\star$;
\item[(ii)]  the restricted gradient $\wt G_{\eps,\delta}(E_\star)$ is nonzero.
\end{itemize}
Let $E(t)\in \C^{n\times n}$ be the solution of \eqref{ode-E} passing through $E_\star$. 
Then the following are equivalent: 
\begin{itemize}
\item[1.] $\displaystyle\frac{ d }{dt} \wt F_{\eps,\delta} \left( E(t) \right)  = 0$. \\
\item[2.] $\dot E = 0$. \\[-3mm]
\item[3.] $E_\star$ is a real multiple of $\wt G_{\eps,\delta}(E_\star)$.
\end{itemize}
\end{theorem}
\medskip

\begin{proof} 3.~implies 2., which implies 1. Further, \eqref{c-s-red} shows that 1. implies 3.
\end{proof}

The following result plays an important role.

\begin{theorem} \label{th:gradprop}
    Under the conditions of Theorem~\ref{thm:stat} we have 
    \begin{equation}
        \text{$\wt G_{\eps,\delta}(E_\star)= \eps G_{\eps,\delta}(E_\star)= - \eps xy^*$},       
    \end{equation}
 where $x$ and $y$ are the left and right eigenvectors to a rightmost eigenvalue of the matrix $A+\eps E_\star+ \delta\eta\,\Pi^\cS E_\star$ with $\eta = 1/\| \Pi^\cS E_\star \|_F$. Hence, $E_\star$ is a real multiple of $xy^*$ and thus of rank 1. 
 Moreover, it follows that $E_\star$ is a stationary point of \eqref{ode-E} if and only if $(E_\star, \eta\,\Pi^\cS E_\star)$  is a stationary point of \eqref{ode-E-ES}.
\end{theorem}
\medskip 

\begin{proof}
    By Theorem~\ref{thm:stat}, $\wt G = \wt G_{\eps,\delta}(E_\star)$ is a real multiple of $E_\star$, and hence (again with $G=- xy^*$ for short)
    \begin{equation} \label{Pi-G}
    \Pi^\cS \wt G =  \eps \Pi^\cS G + \delta \eta \,\Pi^\cS G -
    \delta \eta \, \Re \langle \Pi^\cS G, \eta\Pi^\cS E_\star \rangle \, \eta\Pi^\cS E_\star
    \end{equation}
is a real multiple of $\Pi^\cS E_\star$. This implies that $\Pi^\cS G$ is a real multiple of $\Pi^\cS E_\star$. So
the last two terms in \eqref{Pi-G} cancel, which are the same two terms as in \eqref{red-grad}. This implies 
$\wt G = \eps G$. Moreover, we recall that $E_\star$ is a real multiple of $\wt G$ and thus of $G$, and hence the result follows.
\end{proof}

\section{Rank-1 reduced optimization problem}
\label{sec:rank-one}

\subsection{Rank-1 constrained gradient flow}
\label{subsec:rank1-gradient-flow}
In the differential equation \eqref{ode-E} we project the right-hand side to the tangent space $T_E\cM_1$ of 
the manifold $\cM_1\subset\C^{n\times n}$ of rank-1 matrices at $E\in\cM_1$, 
so that we obtain a differential equation on the rank-1 manifold $\cM_1$. The orthogonal projection $P_E:\C^{n\times n}\to T_E\cM_1$ is known from \cite[Lemma~4.1]{KL07} to be given, for $E=\sigma u v^*\in \cM_1$ with $\sigma>0$ and $u,v\in \C^n$ of unit Euclidean norm, by the expression
$$
P_E(Z) = Z - (I-uu^*) Z (I-vv^*)
\quad\text{ for $Z\in\C^{n, n}$}.
$$
We consider the projected differential equation \eqref{ode-E} on $\cM_1$,
\begin{equation}\label{ode-E-1}
\dot E = -P_E\Bigl( \wt G_{\eps,\delta}(E) - \Re \langle \wt G_{\eps,\delta}(E),E \rangle E \Bigr).
\end{equation}
If the Frobenius norm of $E=E(t)$ is 1, then it is readily checked, using $P_E(E)=E$, that
$
 \Re\, \langle E, \dot E \rangle = 0.
$
Hence, solutions $E(t)$ of  \eqref{ode-E-1} conserve the Frobenius norm 1 for all $t$.

Since we are only interested in solutions of Frobenius norm 1 of \eqref{ode-E-1}, we can write 
$E=uv^*$ with $u$ and $v$ of unit norm (without the extra factor $\sigma$). 
We then have the following differential equations for $u$ and $v$. 

\begin{lemma} [Differential equations for the two vectors]
\label{lem:uv-1}
For an initial value $E(0)=u(0)v(0)^*$ with $u(0)$ and $v(0)$ of unit norm, the solution of \eqref{ode-E-1} is given as
$E(t)=u(t)v(t)^*$, where $u$ and $v$ solve the system of differential equations (for $\wt G=\wt G_{\eps,\delta}(E)$)
\begin{equation}\label{ode-uv}
\begin{array}{rcl}
 \dot u &=& -\tfrac \iu2 \, \Im(u^*\wt Gv)u - (I-uu^*)\wt Gv
\\[1mm]
 \dot v &=& -\tfrac \iu2 \, \Im(v^*\wt G^*u)v - (I-vv^*)\wt G^*u,
\end{array}
\end{equation}
which preserves $\|u(t)\|=\|v(t)\|=1$ for all $t$.
\end{lemma}
\medskip

\begin{proof} The proof is similar to that of Lemma~3.2 of \cite{GLS23}. We include it for the convenience of the reader.
  We introduce the projection $\widehat P_E$ onto the tangent space at $E=uv^*$ of the submanifold of rank-1 matrices of unit Frobenius norm,
 $$
 \widehat P_E(\wt G) = P_E(\wt G-\Re\langle \wt G,E \rangle E)=P_E(\wt G)-\Re\langle \wt G,E \rangle E.
 $$
 We find
 \begin{align*}
 \widehat P_E(\wt G) &= \wt Gvv^* - uu^*\wt Gvv^* + uu^*\wt G -\Re\langle \wt G,uv^* \rangle uv^*
 \\
 &= (I-uu^*)\wt Gvv^* + uu^*\wt G(I-vv^*) + uu^*\wt Gvv^* - \Re(u^*\wt Gv)uv^*
 \\
 &= (I-uu^*)\wt Gvv^* + uu^*\wt G(I-vv^*) + \iu \,\Im(u^*\wt Gv)uv^*
 \\
 &= \Bigl(  \tfrac \iu2 \, \Im(u^*\wt Gv)u + (I-uu^*)\wt Gv\Bigr) v^* + u\Bigl( \tfrac\iu2\, \Im(u^*\wt Gv)v^* + u^*\wt G(I-vv^*) \Bigr).
 \end{align*}
 For $\dot E=\dot u v^* + u {\dot v}^*$ we thus have $\dot E=- \widehat P_E(\wt G)$ if $u$ and $v$ satisfy \eqref{ode-uv}.
 Since then $\Re(u^*\dot u)=0$ and $\Re(v^*\dot v)=0$, the unit norm of $u$ and $v$ is preserved.
\end{proof}

%
%
%
%
%
%
%
%
%

\subsection{Monotone decay and stationary points}
The projected differential equation \eqref{ode-E-1} has the same monotonicity property as the differential equation \eqref{ode-E}.

\begin{theorem}[Monotone decay of the restricted functional] 
\label{thm:monotone-1}
Let $E(t)\in \C^{n\times n}$ of unit Frobenius norm satisfy the differential equation {\rm (\ref{ode-E-1})}.
Assume that the rightmost eigenvalue $\lambda(A+\eps E(t)+\delta\eta(t)\,\Pi^\cS E(t))$ with $\eta(t) = 1/\| \Pi^\cS E(t) \|_F$ is simple and depends continuously on $t$.
Then,
\begin{equation} 
\frac{d}{dt} \wt F_{\eps,\delta} (E(t)) \le  0.
\label{monotone-1}
\end{equation}
\end{theorem}

 \begin{proof} 
From \eqref{ode-E-1} and using $\dot E \in T_E\cM_1$ and $\Re\langle E, \dot E\rangle =0$ we obtain that
$$
\| \dot E \|_F^2 = - \Re\, \bigl\langle P_E\bigl(\wt G-\Re\langle \wt G,E\rangle E\bigr),\dot E \bigr\rangle
= - \Re\, \bigl\langle \wt G-\Re\langle \wt G,E\rangle E,\dot E \bigr\rangle = - \Re\,\langle \wt G, \dot E \rangle.
$$
Lemma~\ref{lem:red-gradient} and \eqref{ode-E-1} thus imply
$$
\frac1{\kappa } \,\frac{d}{dt} \wt F_{\eps,\delta}(E(t))  
= - \bigl\| P_E\bigl(\wt G-\Re\langle \wt G,E\rangle E\bigr)\bigr\|_F^2 ,
$$
which yields \eqref{monotone-1}.
\end{proof}

Comparing the differential equations \eqref{ode-E} and \eqref{ode-E-1} immediately shows that every stationary point of \eqref{ode-E} is also a stationary point of the projected differential equation \eqref{ode-E-1}. The converse is also true for the stationary points $E$ of unit Frobenius norm with ${P_E(G_{\eps,\delta}(E))\ne 0}$. 

\begin{theorem}[Stationary points]
\label{thm:stat-1}
We assume the conditions of Theorem~\ref{thm:stat}.
\begin{enumerate}
\item Every stationary point $E_\star$ of \eqref{ode-E}, which then is a real multiple of $G_{\eps,\delta}(E_\star)$ and satisfies $\wt G_{\eps,\delta}(E_\star)=  G_{\eps,\delta}(E_\star)$,
is also a stationary point of \eqref{ode-E-1}.
\item Conversely, let $E_\star$ be a stationary point of \eqref{ode-E-1} that satisfies 
$\wt G_{\eps,\delta}(E_\star)=  G_{\eps,\delta}(E_\star)$ and $P_E(G_{\eps,\delta}(E))\ne 0$. Then, $E_\star$ is also a stationary point of \eqref{ode-E}.
\end{enumerate}
\end{theorem}
\medskip

\begin{proof} We show that $E$ is a nonzero real multiple of $\wt G = \wt G_{\eps,\delta}(E)$. By Theorem~\ref{thm:stat}, $E$ is then a stationary point of the differential equation \eqref{ode-E}.

For a stationary point $E$ of \eqref{ode-E-1}, the righthand side must vanish, which shows that $P_E(\wt G)$  is a nonzero real multiple of $E$. Hence, in view of $P_E(E)=E$, we can write $G$ as
$$
G=\mu E + W, \quad\text{ where $\mu\ne 0$ is real and $P_E(W)=0$.}
$$
Since $E$ is of rank 1 and of unit Frobenius norm, $E$ can be written as $E=uv^*$ with $\| u \|=\|v\|=1$. We then have
 $$
 W=W-P_E(W)= (I-uu^*)W(I-vv^*).
 $$
 On the other hand, $G=-xy^*$  is also of rank 1. So we have
 $$
- xy^* = \mu uv^* + (I-uu^*)W(I-vv^*).
 $$
 Multiplying from the right with $v$ yields that $x$ is a complex multiple of $u$, and multiplying from the left with $u^*$ yields that $y$ is a complex multiple of $v$. Hence, $G$ is a complex multiple of $E$. Since we already know that $P_E(G)$ is a nonzero real multiple of $P_E(E)=E$, it follows that $G$ is the same real multiple of $E$. By Theorem~\ref{thm:stat}, $E$ is therefore a stationary point of the differential equation \eqref{ode-E}.
\end{proof}

\begin{remark}\label{rem:exceptional}
\rm
We discuss the nondegeneracy condition $P_E(G_{\eps,\delta}(E))\ne 0$. 
Recall that $G=G_{\eps,\delta}(E)=-xy^*$, where $x$ and $y$ are left and right eigenvectors, respectively, to the eigenvalue $\lambda(A+\eps E+\delta E^\cS)$ with $E^\cS=\Pi^\cS E / \| \Pi^\cS E \|_F$. We consider the situation in which we have $P_E(G)= 0$.
For $E=uv^*$, $P_E(G)= 0$ implies $G=(I-uu^*)G(I-vv^*)$, which yields $Gv=0$ and $u^*G=0$ and therefore $y^*v=0$ and $u^*x=0$. So we have $Ey=0$ and $x^*E=0$. This implies that $\lambda$ is already an eigenvalue of $A+\delta E^\cS$ with the same left and right eigenvectors $x,y$ as for $A+\eps E+ \delta E^\cS$, which is a very exceptional situation.
\end{remark}





\section{Numerical integration by a splitting method}
\label{sec:proto-num}
%
We need to integrate numerically the differential equations (\ref{ode-uv}).
The goal is not to follow a particular solution accurately, but to compute a stationary point.
The simplest method is the normalized Euler method, or normalized gradient descent method; 
however, we find that a more efficient method is obtained with a  splitting method.

\subsection{Preparation}
We abbreviate $\wt G = \wt G_{\eps,\delta}$ and  $G = G_{\eps,\delta}$ and rewrite \eqref{ode-uv} (with $\gamma=u^* \wt G v$) as 
\begin{equation}\label{ode-uv2}
\begin{array}{rcllr}
 \dot u &=& \displaystyle {}-\tfrac \iu2 \, \Im(\gamma)u + \gamma u - \wt G v & = & \tfrac \iu2 \, \Im(\gamma)u + \Re(\gamma) u - \wt G v
\\[2mm]
 \dot v &=&  \displaystyle \tfrac \iu2 \, \Im(\gamma)v + \conj{\gamma}\,v - \wt G^*u & = & {}-\tfrac \iu2 \, \Im(\gamma)v + \Re(\gamma)\,v - \wt G^*u,
\end{array}
\end{equation}

Our normalization of the eigenvectors $x(t)$ and $y(t)$ of $A+\eps E(t)+ \delta E^\cS(t)$ 
is that they are of unit norm and with real and positive inner product.
This still leaves us a degree of freedom in the choice of the phase of one of them, say $x(t)$.
We use that to ensure that 
\begin{equation} \label{eq:ass2}
x(t)^* u(t) \in \R^+
\end{equation}
after which we rotate $y$ such that also $x(t)^* y(t) \in \R^+$.

Under the conditions of Theorem~\ref{thm:stat} we have that at a stationary point
$\wt G = - \eps xy^*$,   which implies (with $\alpha=x^*u$ and $\beta=y^*v$)
\[
\left(  \tfrac \iu2 \, \Im(\gamma) + \Re(\gamma) \right) u \propto \beta x, \qquad
\left( -\tfrac \iu2 \, \Im(\gamma) + \Re(\gamma) \right) v \propto 
y 
\]    
where $\propto$ indicates real proportionality. 
This implies
$\tfrac \iu2 \, \Im(\gamma) + \Re(\gamma) \propto \beta$.

We know that at such a stationary point $u v^* \propto x y^*$. 
This implies $\gamma \in \R$ and thus $\beta \in \R$. This suggests the following
splitting of the right-hand side of \eqref{ode-uv2}.

\subsection{Splitting}
The splitting method consists of a first step applied to the differential equations
(with $\gamma=u^* \wt G v$)
\begin{equation}\label{ode-uv-horiz}
\begin{array}{rcl}
\dot u &=&  \Re(\gamma) u - \wt G v  
\\[2mm]
\dot v &=&  \Re(\gamma)\,v - \wt G^*u   
\end{array}
\end{equation}
followed by a step for the differential equations 
\begin{equation}\label{ode-uv-rot} 
\begin{array}{rcl}
\dot u &=&    \displaystyle \tfrac \iu2 \, \Im(\gamma) u 
\\[2mm]
\dot v &=& {}-  \displaystyle \tfrac \iu2 \, \Im(\gamma) v.
\end{array}
\end{equation}

Note that the second differential equation is a rotation of $u$ and $v$. 
In  the case of a real eigenvalue of a real matrix,  the system \eqref{ode-uv-rot} has a vanishing 
right-hand side and can therefore be ignored.



%

Next we get that under the condition on the gradient considered in Theorem \ref{thm:stat-1}, the splitting method preserves stationary points.
\begin{lemma}[Stationary points] \label{lem:stat-split}
Let  $(u,v)$ be a stationary point of the differential equations \eqref{ode-uv} such that 
$E_\star = u v^*$ (which is a stationary point of \eqref{ode-E-1}) satisfies 
$\wt G_{\eps,\delta}(E_\star)=  G_{\eps,\delta}(E_\star)$ and moreover $P_E(G_{\eps,\delta}(E))\ne 0$.
Then $(u,v)$ is a stationary point of the differential equations \eqref{ode-uv-horiz} and 
\eqref{ode-uv-rot}. 

Conversely, let $(u,v)$ be a stationary point of the differential equations \eqref{ode-uv-horiz} and 
\eqref{ode-uv-rot}, such that $E_\star = u v^*$ satisfies $\wt G_{\eps,\delta}(E_\star)=  G_{\eps,\delta}(E_\star)$ and  
$P_E(G_{\eps,\delta}(E))\ne 0$.
Then $(u,v)$ is a stationary point of the differential equations \eqref{ode-uv}.
\end{lemma}

\begin{proof} 
The proof is analogous to that of \cite[Lemma 5.1]{GLS23}.  
\end{proof}

\subsection{Fully discrete splitting algorithm} \label{sec:splitting}

Starting from vectors $u_k,v_k$ of unit norm,
we denote by $x_k$ and $y_k$ the left and right normalized eigenvectors with positive inner product to the target 
eigenvalue $\lambda_{k}$ of $A+\eps E _k+\delta \eta_k \Pi^\cS E_k$, with
$\eta_k = 1/\| \Pi^\cS E_k \|_F$, and set 
\begin{equation}\label{setG-gamma}
G_k = x_k y_k^*, \quad \wt G_k = \eps G_k + \delta\eta_k\, \Pi^\cS G_k - 
\delta \eta_k \,\Re\bigl\langle G_k,\eta_k\,\Pi^\cS E_k\bigr\rangle \,\eta_k\,\Pi^\cS E_k, \quad
\gamma_k = u_k^* \wt G_k v_k. 
\end{equation}
We apply a step of the Euler method with step size $h$  to \eqref{ode-uv-horiz} to obtain
\begin{equation}
\label{eul-horiz}
\begin{array}{rcl}
{\widehat u}(h) &=& u_k + h \left( \Re(\gamma_k)\, u_k - \wt G_k \,v_k \right)
\\[2mm]
{\widehat v}(h) &=& v_k + h \left(\Re(\gamma_k)\, v_k - \wt G_k^* \,u_k \right),
\end{array}
\end{equation}
followed by a normalization to unit norm 
\begin{equation} \label{eq:normal}
\check u(h)=\frac{\widehat u(h)}{\|\widehat u(h)\|},\quad
\check v(h)=\frac{\widehat v(h)}{\|\widehat v(h)\|}.
\end{equation}

Then, as a second step, we integrate the  rotating differential equations \eqref{ode-uv-rot} by setting, 
with $\check\gamma_k = \check u(h)^* \wt G_k \check v(h)$, $\vartheta = \displaystyle \tfrac 12\, \Im(\check\gamma_k)$,
\begin{equation} \label{eq:rotate}
u(h)=\e^{\iu \vartheta h} \, \check u(h), \qquad
v(h)=\e^{{}-\iu \vartheta h} \, \check v (h),
\end{equation} 
which would be solved exactly if $\vartheta$ were constant.

In our numerical experiments, the differential equations for $u$ and $v$ are solved numerically using the proposed 
splitting method possibly coupled with an Armijo-type stepsize selection.

One motivation for choosing this method is that near a {\em non real} stationary point, the motion appears almost 
rotational. 
This algorithm requires in each step one computation of the target eigentriple
of structure-projected rank-$1$ perturbations to the matrix $A$, which can be computed at moderate computational cost for a
large sparse matrix $A$ by a Krylov Schur algorithm \cite{St01}, implemented in the \verb@MATLAB@ function \verb@eigs@.  

\section{Outer iteration}

In the outer iteration we compute $\delta_\eps$, the smallest positive solution of the one-dimensional root-finding problem \eqref{zero-delta}.
This can be solved by a variety of methods, such as bisection. We aim for a locally quadratically convergent Newton-type method, which can be justified under regularity assumptions that appear to be usually satisfied 
If these are not met, we can always resort to bisection. The proposed algorithm in fact uses a combined 
Newton / bisection approach.

For a fixed $\eps>0$ we define
\begin{equation} \label{Eoptim}
E(\delta) = E_{\eps,\delta} = \arg \min\limits_{E \in \C^{n,n}} \widetilde F_{\eps,\delta}\left( E \right).
\end{equation}

\begin{assumption} \label{ass:E-delta}
For $\delta$ close to $\delta_\eps$, 
we assume the following for $E(\delta)$:
\begin{itemize}
\item[(i) ]   The rightmost eigenvalue $\lambda(\delta)=\lambda\bigl( A+\eps E(\delta) + \delta\, \eta(\delta)\,\Pi^\cS E(\delta) \bigr)$ with $\eta(\delta)=1/\|\Pi^\cS E(\delta)\|_F$ is a simple eigenvalue.
\item [(ii) ] The map $\delta \mapsto E(\delta)$ is continuously differentiable.
\end{itemize}
\end{assumption}
Under this assumption, the branch of eigenvalues $\lambda(\delta)$ and its corresponding eigenvectors $x(\delta), y(\delta)$ normalized such that they are of unit norm and with real and positive inner product are also continuously differentiable functions of $\delta$ in a left neighborhood of 
$\delta_\eps$. 
We denote  the eigenvalue condition number by
\[
\kappa(\delta) = \frac 1 {x(\delta)^*y(\delta)}>0.
\]
The following result gives us an explicit and easily computable expression for the derivative of the function
$\ophi_\eps(\delta) = \widetilde F_{\eps,\delta}\left( E(\delta) \right) = -\Re \,\lambda( \delta)  $ 
with respect to $\delta$.

\begin{theorem}[Derivative for the Newton iteration] 
\label{thm:phi-derivative}
Under Assumption~\ref{ass:E-delta}, the function $\ophi(\delta)$ is continuously differentiable in a neighborhood of $\delta$ 
and its derivative is given as
\begin{equation} \label{derdelta}
\ophi_\eps'(\delta)  = 
-\kappa(\delta) \, \| \Pi^\cS (x(\delta)y(\delta)^*) \|_F .
\end{equation}
\end{theorem}
\begin{proof}
We use Lemma~\ref{lem:gradient} and proceed similarly to the proof of Lemma \ref{lem:red-gradient}.
Indicating by $'$ differentiation w.r.t.~$\delta$ 
and noting that, with $\eta(\delta) = 1/\| \Pi^\cS E(\delta) \|_F$,
\[
\frac d{d\delta}\left(\delta\, \eta(\delta) \,\Pi^\cS E(\delta) \right) = \eta(\delta)\, \Pi^\cS E(\delta) + \delta \,\eta'(\delta) \,\Pi^\cS E(\delta) + \delta \,\eta(\delta)\, \Pi^\cS E'(\delta),
\]
we obtain with the gradient $G(\delta)=- x(\delta)y(\delta)^*$ (cf.~Lemma~\ref{lem:gradient}) and with $\wt G(\delta)$ defined as in \eqref{red-grad} with $E(\delta)$ and $G(\delta)$,
\begin{equation} \label{derivdelta}
\frac{1}{\kappa(\delta)} \,\frac{d}{d \delta} \widetilde F_{\eps,\delta}(E(\delta)) = 
\Re \bigl\langle  G(\delta),  \eta(\delta){\Pi^\cS E(\delta)}  \bigr\rangle
+ \Re \bigl\langle  \widetilde G(\delta),  E'(\delta) \bigr\rangle.
\end{equation}
We know by Theorem \ref{thm:stat-red} that in the stationary point $E(\delta)$, there exists a real $\mu(\delta)$ such that
$$
E(\delta) = \mu(\delta) \widetilde G(\delta).
$$
Since $\| E(\delta) \|_F=1$ for all $\delta$, we find  $1= |\mu(\delta)| \, \| \widetilde G(\delta)\|_F $ (in particular $\mu(\delta)\ne 0$) and
\[
0= \frac12 \frac d{d\delta} \| E(\delta) \|^2 = \Re \langle E(\delta),E'(\delta) \rangle = \mu(\delta)\, \Re \langle \widetilde G(\delta), E'(\delta) \rangle,
\]
so that the last term in \eqref{derivdelta} vanishes.
Using that $\widetilde G(\delta) = \eps G(\delta)$ by Theorem \ref{th:gradprop}, we obtain
\[
\frac{1}{\kappa(\delta)} \ophi_\eps'(\delta) = \Re \biggl\langle  G(\delta),  \pm\frac{\Pi^\cS G(\delta)}{\| \Pi^\cS G(\delta) \|} \biggr\rangle = \Re \biggl\langle  \Pi^\cS G(\delta),  \pm\frac{\Pi^\cS G(\delta)}{\| \Pi^\cS G(\delta) \|}  \biggr\rangle  =  \pm\| \Pi^\cS G(\delta) \|_F. 
\]
Since $\phi_\eps$ is monotonically decreasing, the correct sign is the minus sign. This yields the stated result.
\end{proof}

For the converse problem \eqref{zero-eps} of finding, for a given $\delta>0$, the zero $\eps_\delta$ of  $\psi_\delta(\eps)=\wt F_{\eps,\delta}(E_{\eps,\delta})$, we obtain in the same way
\begin{equation}\label{dereps}
\psi_\delta'(\eps) = - \kappa(\eps) \,\| G(\eps) \|_F = - \kappa(\eps) <0.
\end{equation}

\begin{algorithm}
\DontPrintSemicolon
\KwData{Matrix $ A$, matrix type (real/complex, structured), $\eps>0$,\; 
${\rm tol}_0$ (initial tolerance), $k_{\max}$ (max number of iterations)\; 
$\delta_{\rm lb}$ and $\delta_{\rm ub}$ (starting values for the lower and upper 
bounds for $\delta_\eps$)} 
\KwResult{$\delta_\eps$ (computed value / upper bound of the $\eps$-stability radius)}
\Begin{
\nl Set $\lambda(0)$ rightmost eigenvalue of $ A$, $x(0)$ and $y(0)$ the corresponding left and right
eigenvectors of unit norm with $x(0)^*y(0)>0$.\; 
\nl Initialize $\delta_0$ and $E(\delta_0)$ according to the setting. \; 
    Set $k=0$.\;
\nl Initialize lower and upper bounds: $\delta_{\rm lb}=0$, $\delta_{\rm ub}=+\infty$.\; 
\While{$|\Re\,\lambda(\delta_k)-\Re\,\lambda(\delta_{k-1})| < {\rm tol}_k$ {\bf and} $k<k_{\max}$}{
\nl {\bf Inner iteration:} Compute $E(\delta_k)=u(\delta_k)v(\delta_k)^*$ by integrating the rank-1 matrix differential equation \eqref{ode-uv} 
with initial datum $E(\delta_{k-1})$ into a stationary point. This also yields a rightmost eigenvalue $\lambda(\delta_k)$ of
$A+\eps E(\delta_k)+ \delta_k \Pi^\cS E(\delta_k)/\|\Pi^\cS E(\delta_k)\|_F$ and its left and right eigen\-vectors
$x(\delta_k)$ and $y(\delta_k)$ of unit norm with $x(\delta_k)y(\delta_k)^*>0$.\;
\nl Update upper and lower bounds $\delta_{\rm lb}$, $\delta_{\rm ub}$:\; 
\nl \eIf{$\Re\,\lambda(\delta_k) > 0$} 
{Set $ \delta_{\rm ub} = \min(\delta_{\rm ub},\delta_k)$.} 
{Set $ \delta_{\rm lb} = \max(\delta_{\rm lb},\delta_k)$.}
\nl Compute $\delta_{k+1} = \delta_{k} -
\dfrac{x(\delta_k)^* y(\delta_k)}{\| \Pi^\cS (x(\delta)y(\delta)^*) \|_F}\, \Re\,\lambda(\delta_k)$.\;
\nl Set $k=k+1$.\; 
\nl \If{$\delta_{k} \not\in [\delta_{\rm lb},\delta_{\rm ub}]$}
{Set $\delta_{k} = (\delta_{\rm lb} + \delta_{\rm ub})/2$.}
{Set ${\rm tol}_k = \max \{ 10^{-2}\,{\rm tol}_{k-1}, 10^{-8} \}$.}\;
}
\nl \eIf{$k \le k_{\max}$}
{Set $\delta_\eps = \delta_k$.}{Print {\em max number of iterations reached.}}
}
\caption{Outer iteration: Newton / bisection method}
\label{alg_SR} 
\end{algorithm}

Algorithm \ref{alg_SR} implements a hybrid Newton / bisection method that maintains
an interval known to contain the root, bisecting when the Newton step is outside
the interval $[ \delta_{\rm lb}, \delta_{\rm ub} ]$.

Due to the possible convergence of the inner method to a local instead of global minimum, the final value $\delta_\eps$
computed by Algorithm \ref{alg_SR} might be larger than the minimal one. The computed value of $\delta_\eps$ is thus an upper bound. Computing several trajectories with different starting values reduces the risk of getting caught in a nonoptimal local minimum. Our numerical experience indicates, however, that the global minimum is usually found by a single discrete trajectory starting with
$E(0)=xy^*$, where $x$ and $y$ are left and right eigenvectors of the unperturbed matrix $A$, of unit norm and with positive inner product.


\section{Numerical experiments}

We first apply our method on a small-size Grcar matrix 
and then on a few sparse examples, considered e.g. in \cite{R15}, the Tolosa matrix of dimension 
$n=4000$ (see e.g. \cite{BPR97,W02}). 
Finally we consider the Tubular matrix of dimension $n=1000$ \cite{BPR97}, on which we fix $\delta$ and compute the smallest
common bound for the resolvent norms of the perturbed matrices with sparsity-structured perturbations of Frobenius norm at most $\delta$.
All considered examples are characterized by complex conjugate rightmost eigenvalues, except for the last one.

In our experiments we have obtained the extremizers always with the positive sign $+\delta$ in \eqref{eq:redF}.

\subsection*{Numerical considerations}

Based on our experience, for very sparse matrices of dimension $n \gg 1$, in terms of CPU time, 
it is convenient to use the full problem instead of its rank-1 projection and possibly exploit 
the sparse plus rank-1 structure of the matrices in the eigenvalue computation, through Krylov 
subspace methods (as done by \verb@eigs@). 
The number of steps of the standard Euler integrator in fact turns out to be smaller than the 
one for the splitting rank-1 method in the experiments we made.

For dense structures and for sparse structures with a number of nonzero entries equal to $c n$ 
(with $c$ significantly larger than $1$) it is convenient to use the rank-$1$ system solved by 
the splitting integrator.

In our implementation, the computation of the eigenvalues in very sparse examples is   
achieved by the \verb@MATLAB@ routine \verb@eigs@  \cite{St01,LeSY98} with a default choice of $20$ 
Ritz values. The convergence of \verb@eigs@ is not always guaranteed and in a numerical implementation 
its use has to be considered carefully.

\subsection{A small illustrative example}
\label{sec:small}

We take the matrix $A = -{\rm Grcar}(10) - I$ ($I$ stands for the identity matrix) 
from the Eigtool demo,
that is
\[
A = 
\left( \begin{array}{rrrrrrrr}
  -2  & -1  & -1  & -1  &  0  &  0  &  0  & 0  \\ 
   1  & -2  & -1  & -1  & -1  &  0  &  0  & 0  \\
   0  &  1  & -2  & -1  & -1  & -1  &  0  & 0  \\
   0  &  0  &  1  & -2  & -1  & -1  & -1  & 0  \\
   0  &  0  &  0  & \ddots & \ddots & \ddots & \ddots & \ddots	
\end{array}
\right).
\]
\normalsize
We fix $\eps = 0.5$ (in which case the $\eps$-pseudospectral abscissa of $A$ is given by 
$\alpha_\eps(A) = {}-3.890782704837603 \cdot 10^{-1}$) and first choose $\cS$ as the space 
of real matrices with the same sparsity pattern of $A$. 

The rightmost eigenvalues of $A$ are complex conjugate,
\[
\lambda = -1.197971039973676  \pm 2.129259562786844 \iu.
\] 
and the stability radius turns out to be given by
$
\oeps =  8.39282612 \cdot 10^{-1},
$
which provides  a lower bound for the structured $\eps$-stability radius $\delta_\eps^\cS$: 
$$
\delta_{\eps}^\cS \ge  \oeps - \eps = 3.39282612 \cdot 10^{-1},
$$
which is the unstructured $\eps$-stability radius.

\begin{table}[ht]
\begin{center}
\begin{tabular}{|l|l|r|l|}\hline
  $k$ & $\delta_k$ & $\Re\,\lambda_k$ & $\#$ steps \\
 \hline
\rule{0pt}{9pt}
\!\!\!\! 
  $1$         & $0$ & $ {}-3.8907827  \cdot 10^{-1}$   & $110$  \\
	$2$         & $8.5881368 \cdot 10^{-1}$ & $ 3.0135918  \cdot 10^{-3}$   & $126$  \\
	$3$         & $8.5228455 \cdot 10^{-1}$ & $ 3.3695994  \cdot 10^{-7}$   & $94$  \\
	$4$         & $8.5228382 \cdot 10^{-1}$ & $ 3.4902636  \cdot 10^{-15}$   & $5$  \\
 \hline
\end{tabular}
\vspace{2mm}
\caption{Computed values of $\delta_k$, $\phi_\eps(\delta_k)$ and number of steps of the splitting method 
using Algorithm {\rm \ref{alg_SR}} for the shifted Grcar matrix of dimension $10$.
\label{tab:exG10S}}
\end{center}
\end{table}


\begin{figure}[ht]
\begin{center}
\includegraphics[height=8cm,width=8.6cm]{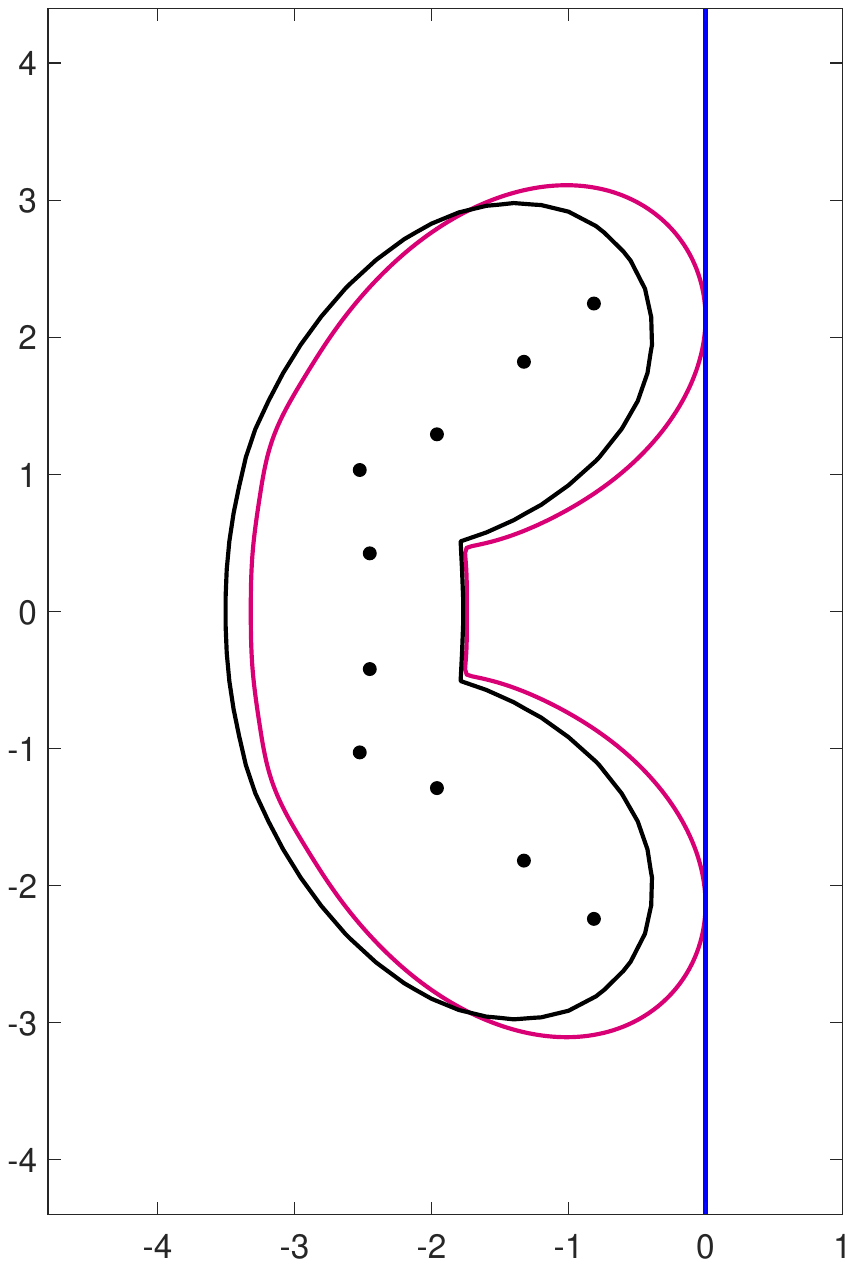} 
\vspace{-1cm}
\caption{The $\eps$-pseudospectra of the matrix $A$ (in black) of subsection {\rm \ref{sec:small}} and $A + \Delta_\eps^\cS$
(in red) for the value $\eps=0.5$ and $\cS$ the space of real sparse matrices with the sparsity pattern of $A$.  The vertical
line shows the imaginary axis. The black bullets indicate the eigenvalues of~$A$. 
}
\label{fig:G10}
\end{center}
\end{figure}

We report the results obtained applying Algorithm \ref{alg_SR} coupled with the splitting method
described in subsection \ref{sec:splitting} in Table \ref{tab:exG10S}.
The legend is the following: $k$ addresses 
the number of the outer iteration, $\delta_k$ is the norm of the structured perturbation,
$\phi_\eps(\delta_k)$ is the largest  $\eps$-pseudospectral abscissa and $\# $ steps is the number
of steps executed by the splitting integrator, which identify the number of calls to either the routine 
\verb@eig@ (for dense problems) or \verb@eigs@ (for sparse problems), which yields the most expensive
part of the whole algorithm. 

We find $\delta_\eps^\cS = 8.5228382298260 \cdot 10^{-1}$. The $\eps$-pseudospectra of $A$ and 
$A + \Delta_\eps^\cS$ with $\Delta_\eps^\cS=\delta_\eps^\cS \Pi^\cS E(\delta_\eps^\cS)/\|\Pi^\cS E(\delta_\eps^\cS)\|_F$ are shown in Figure \ref{fig:G10}. Applying the criss-cross
algorithm of \cite{BuLeOv03} confirms that the computed $\delta_\eps^\cS$ and the computed matrix $\Delta_\eps^\cS$ provides
the global maximum of the optimization problem to find $\delta>0$ and $\Delta\in\cS$ such that
\eqref{robust-res-bound} is satisfied, viz.,
\[
  \max_{\Delta\in\cS, \| \Delta\|_F \le \delta}\max_{\,\Re\lambda \ge 0\,} \| (A+\Delta-\lambda I)^{-1} \|_2 = \frac 1\eps.
\]

Next we consider the dual problem, fixing $\delta = 8.5228382298260 \cdot 10^{-1}$ (the $\delta$-bound we have just computed) and look for $\eps_\delta^\cS$ such that the $\delta$-robust resolvent bound \eqref{robust-res-bound} on the right half-plane is satisfied with $\eps=\eps_\delta^\cS$.
This should provide the value $\eps=0.5$. The results are indeed striking:
\begin{table}[ht]
\begin{center}
\begin{tabular}{|l|l|r|l|}\hline
  $k$ & $\eps_k$ & $\Re\,\lambda_k$ & $\#$ steps \\
 \hline
\rule{0pt}{9pt}
\!\!\!\! 
        $1$         & $                   10^{-2}$ & ${}-7.5342048842 \cdot 10^{-1}$    & $657$  \\
	$2$         & $3.5214566707 \cdot 10^{-1}$ & ${}-1.9016325803 \cdot 10^{-1}$    & $170$  \\
	$3$         & $4.9382789893 \cdot 10^{-1}$ & ${}-7.6628326477 \cdot 10^{-3}$    & $119$  \\
	$4$         & $4.9999196551 \cdot 10^{-1}$ & ${}-9.9621473717 \cdot 10^{-6}$    & $91$  \\
	$5$         & $4.9999999999 \cdot 10^{-1}$ & ${}-1.6626651445 \cdot 10^{-11}$   & $51$  \\
 	$6$         & $5 \cdot 10^{-1}$            & ${}-2.0816681712 \cdot 10^{-16}$   & $2$  \\
 \hline
\end{tabular}
\vspace{2mm}
\caption{Computed values of $\eps_k$, $\Re\,\lambda_k$ and number of steps of the numerical integrator
using the analog of Algorithm {\rm \ref{alg_SR}} for the shifted Grcar matrix.
\label{tab:exG10E}}
\end{center}
\end{table}

Finally -- for the same problem -- we choose $\cS$ as the set of Toeplitz real matrices with the same band of $A$, for which
we obtain the results in Table \ref{tab:exG10T}.
The extremal structured perturbation (of norm $9.043542933808467 \cdot 10^{-1}$) is given by
\[
\Delta = 10^{-2}\, {\rm T}_{1,3} \left( 3.295829030,   7.282237246,     2.619428085,   -4.166704981,    -4.668125451 \right)
\]
where ${\rm T}_{1,3} (a,b,c,d,e)$ denotes the Toeplitz matrix with diagonals from $-1$ to $3$ with diagonal entries
ordered according to the increasing index of the diagonal ($a$ in the lower sub-diagonal, $b$ in the main diagonal, $c, d$ and $e$ in the upper
super-diagonals).

\begin{table}[ht]
\begin{center}
\begin{tabular}{|l|l|r|l|}\hline
  $k$ & $\delta_k$ & $\Re\,\lambda_k$ & $\#$ steps \\
 \hline
\rule{0pt}{9pt}
\!\!\!\! 
  $1$         & $0$ & $ {}-3.8907827  \cdot 10^{-1}$   & $110$  \\
	$2$         & $9.0397503545 \cdot 10^{-1}$ & $ -1.6314020426  \cdot 10^{-4}$   & $125$  \\
	$3$         & $9.0435429340 \cdot 10^{-1}$ & $  6.2862705492  \cdot 10^{-12}$   & $67$  \\
	$4$         & $9.0435429338 \cdot 10^{-1}$ & $ -3.5527136788  \cdot 10^{-15}$   & $2$  \\
 \hline
\end{tabular}
\vspace{2mm}
\caption{Computed values of $\delta_k$, $\Re\,\lambda_k$ and number of steps of the splitting method 
using Algorithm~{\rm \ref{alg_SR}} for the Grcar matrix of dimension $10$, with $\cS$ the space of real 
sparse matrices with the sparsity pattern of $A$.
\label{tab:exG10T}}
\end{center}
\end{table}

\subsection{The Tolosa matrix}
\label{sec:Tolosa}

We continue by considering the Tolosa matrix of dimension 
and $n=4000$ (see e.g. \cite{BPR97,W02}). We fix $\eps=10^{-3}$.
This is a real sparse Hurwitz matrix with righmost complex conjugate eigenvalues.

In Table \ref{tab:exT4000} we report the results of Algorithm \ref{alg_SR} for the Tolosa matrix of dimension $n=4000$.
The quadratic convergence in the outer iteration is reached rapidly and the numerical integrator converges to the
stationary point in a small number of steps.

\begin{table}[ht]
\begin{center}
\begin{tabular}{|l|l|r|l|}\hline
  $k$ & $\delta_k$ & $\Re\,\lambda_k$ & $\#$ steps \\
 \hline
\rule{0pt}{9pt}
\!\!\!\! 
  $1$         & $                0$          & ${}-7.7992086890 \cdot 10^{-2}$    & $2$  \\
	$2$         & $1.5598421556 \cdot 10^{-1}$ & $   2.4138894085 \cdot 10^{-4}$    & $30$  \\
	$3$         & $1.5550522950 \cdot 10^{-1}$ & $   1.1407707057 \cdot 10^{-6}$    & $5$  \\
	$4$         & $1.5550296584 \cdot 10^{-1}$ & $   5.3690606519 \cdot 10^{-9}$    & $2$  \\
	$5$         & $1.5550295518 \cdot 10^{-1}$ & $   2.5638480761 \cdot 10^{-11}$   & $3$  \\
	$6$         & $1.5550295513 \cdot 10^{-1}$ & ${}-1.7186173186 \cdot 10^{-13}$   & $2$  \\
 \hline
\end{tabular}
\vspace{2mm}
\caption{Computed values of $\delta_k$, $\Re\,\lambda_k$ and number of steps of the splitting method 
using Algorithm {\rm \ref{alg_SR}} for the Tolosa-$4000$ matrix.
\label{tab:exT4000}}
\end{center}
\end{table}

\subsection{The Tubular matrix}
\label{sec:Tubular} 

We conclude by considering the Tubular matrix of dimension $n=1000$ (see \cite{BPR97}) and let $\cS$ be the
linear space of real matrices with the same sparsity pattern as $A$. 
In this case, differently from previous examples, we fix $\delta=0.1$ and look for $\eps_\delta^\cS$ 
such that the $\delta$-robust resolvent bound \eqref{robust-res-bound} on the right half-plane is satisfied with $\eps=\eps_\delta^\cS$. This is computed with a variant of Algorithm 1, which yields the resolvent bound
$$
\displaystyle\frac{1}{\eps_\delta^\cS} = \frac{1}{1.12242717731 \cdot 10^{-1}} = 8.909263961258253.
$$
\medskip

\begin{table}[ht]
\begin{center}
\begin{tabular}{|l|l|r|l|}\hline
  $k$ & $\eps_k$ & $\Re\,\lambda_k$ & $\#$ steps \\
 \hline
\rule{0pt}{9pt}
\!\!\!\! 
  $1$         & $                   10^{-2}$ & ${}-8.2035030666 \cdot 10^{-1}$    & $36$  \\
	$2$         & $1.0205194415 \cdot 10^{-1}$ & ${}-6.9429177622 \cdot 10^{-2}$    & $202$  \\
	$3$         & $1.1218404803 \cdot 10^{-1}$ & ${}-3.9743783078 \cdot 10^{-4}$    & $161$  \\
	$4$         & $1.1224271349 \cdot 10^{-1}$ & ${}-2.8689715048 \cdot 10^{-8}$    & $6$  \\
	$5$         & $1.1224271772 \cdot 10^{-1}$ & ${}-6.1345151181 \cdot 10^{-11}$   & $6$  \\
 \hline
\end{tabular}
\vspace{2mm}
\caption{Computed values of $\eps_k$, $\Re\,\lambda_k$ and number of steps of the splitting method 
using the analog of Algorithm {\rm \ref{alg_SR}} for the Tubular-$1000$ matrix.
\label{tab:exU1000}}
\end{center}
\end{table}

As final test we consider the dual problem, i.e. we fix $\eps = \eps_\delta^\cS = 0.112242717885079$ 
and look for $\delta_\eps^\cS$, which is obtained applying Algorithm \ref{alg_SR}. 
The deviation from the exact value $\delta = 0.1$ is presumably due to inaccuracy in the eigenvalue 
computation in the inner iteration.

\begin{table}[ht]
\begin{center}
\begin{tabular}{|l|l|r|l|}\hline
  $k$ & $\delta_k$ & $\Re\,\lambda_k$ & $\#$ steps \\
 \hline
\rule{0pt}{9pt}
\!\!\!\! 
        $1$         & $                0$          & ${}-1.2071026698 \cdot 10^{-2}$    & $26$  \\
	$2$         & $1.0133473305 \cdot 10^{-1}$ & $   1.6334901455 \cdot 10^{-4}$    & $119$  \\
	$3$         & $9.9999834364 \cdot 10^{-2}$ & $   5.7277831811 \cdot 10^{-8}$    & $21$  \\
	$4$         & $9.9999833473 \cdot 10^{-2}$ & $   1.0900480518 \cdot 10^{-10}$   & $4$  \\
	$5$         & $9.9999833956 \cdot 10^{-2}$ & ${}-5.9063642865 \cdot 10^{-11}$   & $2$  \\
 \hline
\end{tabular}
\vspace{2mm}
\caption{Computed values of $\delta_k$, $\Re\,\lambda_k$ and number of steps of the splitting method 
using Algorithm {\rm \ref{alg_SR}} for the Tubular-$1000$ matrix with $\eps = \eps_\delta^\cS$.
\label{tab:exD1000}}
\end{center}
\end{table}

\subsection*{Acknowledgments}



Nicola Guglielmi acknowledges that his research was supported by funds from the Italian MUR (Ministero 
dell'Universit\'a e della Ricerca) within the PRIN 2021 Project Advanced numerical methods for time dependent 
parametric partial differential equations with applications and the Pro3 Project Calcolo scientifico per 
le scienze naturali, sociali e applicazioni: sviluppo metodologico e tecnologico. 
Nicola Guglielmi is affiliated to the Italian INdAM-GNCS (Gruppo Nazionale di Calcolo Scientifico).

Christian Lubich acknowledges the hospitality of GSSI in L'Aquila during a visit in June 2023, where this research originated.


\begin{thebibliography}{10}








\bibitem{BuLeOv03}
J.~V.~Burke, A.~S.~Lewis, and M.~L.~Overton.
Robust stability and a criss-cross algorithm for pseudospectra.
IMA J. Numer. Anal., 23(3): 359--375, 2003.

\bibitem{BPR97}
R.~F.~Boisvert, R.~Pozo, K.~Remington, R.~F.~Barrett, and J.~J.~Dongarra.
Matrix Market: A Web Resource for Test Matrix Collections.
Chapman \& Hall
\newblock http://math.nist.gov/MatrixMarket/, 1997.

\bibitem{GreLO20}
A.~Greenbaum,  R.-C.~Li, and M.~L.~Overton.
First-order perturbation theory for eigenvalues and eigenvectors.
SIAM Rev., 62(2): 463--482, 2020. 


\bibitem{GL11}
N.~Guglielmi and C.~Lubich.
Differential equations for roaming pseudospectra: paths to extremal
points and boundary tracking.
SIAM J. Numer. Anal., 49: 1194--1209, 2011.




\bibitem{GLS23}
N.~Guglielmi, C.~Lubich, and S.~Sicilia.
Rank-$1$ matrix differential equations for structured eigenvalue optimization. 
SIAM J. Numer. Anal., 61: 1737--1762, 2023. 


 \bibitem{GO11}
 N.~Guglielmi and M.~L. Overton.
 Fast algorithms for the approximation of the pseudospectral abscissa
 and pseudospectral radius of a matrix.
 SIAM J. Matrix Anal. Appl., 32(4): 1166--1192, 2011.


\bibitem{HinP05}
D.~Hinrichsen and A.~J.~Pritchard.
Mathematical systems theory I: modelling, state space analysis, stability and robustness.
Springer, Berlin, 2005.






\bibitem{KL07}
O.~Koch and C.~Lubich, Dynamical low-rank approximation.
SIAM J. Matrix Anal. Appl., 29(2): 434--454, 2007.





\bibitem{LeSY98}
R.~B.~Lehoucq, D.~C.~Sorensen, and C.~Yang, ARPACK Users' Guide: Solution of Large-Scale Eigenvalue Problems with Implicitly Restarted Arnoldi Methods, SIAM Publications, Philadelphia, 1998. 







\bibitem{R15}
M.~Rostami.
New algorithms for computing the real structured
pseudospectral abscissa and the real stability radius of
large and sparse matrices.
SIAM J. Sci. Comp., 37(5):447--471, 2015.




\bibitem{St01}
G.~W.~Stewart, 
A {K}rylov-{S}chur algorithm for large eigenproblems,
SIAM J. Matrix Anal. Appl., 23(3), 601--614, 2001/02.





\bibitem{TrE05}
L.~N.~Trefethen, M.~Embree, Spectra and pseudospectra. 
The behavior of nonnormal matrices and operators. Princeton University Press, Princeton, NJ, 2005.


\bibitem{W02}
T.~G.~Wright.
Eigtool: a graphical tool for nonsymmetric eigenproblems.
{\em Oxford University Computing Laboratory,
http://www.comlab.ox.ac.uk/pseudospectra/eigtool/}, 2002.






\end{thebibliography}
\end{document}